\input amstex
\magnification=\magstep1
\baselineskip=13pt
\documentstyle{amsppt}
\vsize=8.7truein
\NoRunningHeads
\def\per{\operatorname{per}}
\def\haf{\operatorname{haf}}
\def\xx{\bold{x}}
\topmatter
\title Integration and Optimization of Multivariate Polynomials by Restriction onto
a Random Subspace   \endtitle
\author Alexander Barvinok \endauthor
\address Department of Mathematics, University of Michigan, Ann Arbor,
MI 48109-1109, USA \endaddress
\email barvinok$\@$umich.edu  \endemail
\date February 2005 \enddate
\thanks This research was partially supported by NSF Grant DMS 0400617.
\endthanks
\abstract We consider the problem of efficient integration of an $n$-variate polynomial 
with respect to the Gaussian measure in ${\Bbb R}^n$ and related problems of
complex integration and optimization of a polynomial on the unit sphere. 
We identify a class of $n$-variate polynomials $f$ for which the integral of any positive
integer power $f^p$ over the whole space
is well-approximated by a properly scaled integral over a random subspace of dimension $O(\log n)$.
Consequently, the maximum of $f$ on the unit sphere is well-approximated by a properly scaled maximum 
on the unit sphere in a random subspace of dimension $O(\log n)$. We discuss connections
with problems of combinatorial counting and applications to efficient approximation of 
a hafnian of a positive matrix.
\endabstract
\endtopmatter
\document
\head 1. Introduction \endhead

We consider the problem of efficient integration of multivariate polynomials with respect to the Gaussian measure in ${\Bbb R}^n$.

Let us assume that the real $n$-variate homogeneous polynomial $f$ of degree $m$ is given to us 
by some ``black box'', which inputs an $n$-vector $x=(\xi_1, \ldots, \xi_n)$ and outputs 
the value of $f(x)$. We want to compute or estimate the integral
$$\int_{{\Bbb R}^n} f \ d \mu_n,$$
where $\mu_n$ is the standard Gaussian measure with the density 
$$(2 \pi)^{-n/2} e^{-\|x\|^2/2}, \quad \text{where} \quad 
\|x\|=\xi_1^2 + \ldots + \xi_n^2 \quad \text{for} \quad x=(\xi_1, \ldots, \xi_n).$$ 
If $m$ is odd then the integral is 0, so the interesting case is that of 
an even degree $m$.

An equivalent problem is to integrate $f$ over the unit sphere ${\Bbb S}^{n-1} \subset {\Bbb R}^n$.
Assuming that $m=2k$ is even, we have 
$$\int_{{\Bbb S}^{n-1}} f(x) \ dx ={\Gamma(n/2) \over 2^k \Gamma(n/2+k) }
 \int_{{\Bbb R}^n} f \ d\mu_n,$$
where $dx$ is the rotation invariant Haar probability measure on ${\Bbb S}^{n-1}$. 
This and related formulas for integrals of polynomials over the unit sphere and over 
the Gaussian measure on ${\Bbb R}^n$ can be found, for example, in \cite{B02b}.

The most straightforward and the most general approach to integration 
is to employ the Monte Carlo method, that is, to sample $N$ random points $x_i \in {\Bbb S}^{n-1}$
and approximate the integral by the sample mean:
$$\int_{{\Bbb S}^{n-1}} f(x) \ dx \approx {1 \over N} \sum_{i=1}^N f(x_i).$$
Although one can show that for a ``typical'' polynomial the Monte Carlo method works reasonably 
well, there are simple examples of polynomials where one would require to sample exponentially many 
points to get reasonably close to the integral.
\example{(1.1) Example} Suppose that $f(x)=\xi_1^{2k}$ for $x=(\xi_1, \ldots, \xi_n)$. Then 
$$\int_{{\Bbb S}^{n-1}}f(x) \ dx= {\Gamma(n/2) \Gamma(1/2 + k) \over \sqrt{\pi} \Gamma(n/2+k)}.$$
If we choose $k \sim n/2$ then the integral is of the order of $2^{-n}$ for large $n$. 

On the other hand, if we sample $N$ random points $x_i$ on the unit sphere ${\Bbb S}^{n-1}$, then 
with high probability we will have $|\xi_1| = O(\sqrt{\ln N / n})$ for the first coordinate
$\xi_1$ of every sampled
point, cf., for example, Section 2 of \cite{MS86}. Thus to approximate the integral within a factor $c^n$ for some absolute constant $c$, the 
number $N$ of samples should be exponentially large in $n$.
\endexample

The reason why the Monte Carlo method doesn't work well on the above example is clear:
the polynomial $f(x)=\xi_1^{2k}$ acquires some large values for an exponentially small fraction 
of $x \in {\Bbb S}^{n-1}$ but those values significantly contribute to the integral. In other words,
the Monte Carlo method wouldn't work well if the graph of the polynomial looks ``needle-like''.
In this paper, we suggest a method tailored specifically for such needle-like polynomials. 

The following defines the class of ``needle-like'' or ``focused''
polynomials we deal with.  
\definition{(1.2) Definitions} Let 
$$\langle x, y \rangle =\xi_1 \eta_1 + \ldots + \xi_n \eta_n \quad 
\text{for} \quad x=(\xi_1, \ldots, \xi_n) \quad \text{and} \quad y=(\eta_1, \ldots, \eta_n)$$
be the standard scalar product in ${\Bbb R}^n$.

Let us fix a number $0 <  \delta \leq  1$ and a positive integer $N$. We say 
that a homogeneous polynomial $f: {\Bbb R}^n \longrightarrow {\Bbb R}$ 
of degree $m$ is $(\delta, N)$-{\it focused} if 
there exist $N$ non-zero vectors $c_1, \ldots, c_N \in {\Bbb R}^n$ such that 
\medskip
$\bullet$ for every pair $(i,j)$ the cosine of the angle between $c_i$ and $c_j$ is at least $\delta$;
\medskip
$\bullet$ the polynomial $f$ can be written as a non-negative linear combination
$$f(x)=\sum_I \alpha_I \prod_{i \in I} \langle c_i, x \rangle, $$
where the sum is taken over subsets $I \subset \{1, \ldots, N\}$ of cardinality 
$|I|=m$ and $\alpha_I \geq 0$.
\enddefinition
Our first result is that the value of the integral of a focused polynomial over a random lower-dimensional subspace allows one to predict the value of the integral over the whole space. 

For a $k$-dimensional subspace $L \subset {\Bbb R}^n$, let $\mu_k$ be the Gaussian 
measure concentrated on $L$ with the density 
$(2 \pi)^{-k/2} \exp\bigl\{-\|x\|^2/2\bigr\}$
for $x \in L$.
We pick a $k$-dimensional subspace at random with respect to the Haar probability measure on
the Grassmannian $G_k({\Bbb R}^n)$ and consider the integral
$$\int_{L} f \ d\mu_k.$$ 
We claim that as long as $k \sim \log N$, the properly scaled integral over $L$ approximates 
the integral over ${\Bbb R}^n$ within a factor of $(1-\epsilon)^{m/2}$.
\proclaim{(1.3) Theorem} There exists an absolute constant $\gamma>0$ with 
the following property.

For any $\delta >0$, for any positive integer $N$, for any  
$(\delta, N)$-focused polynomial \break $f: {\Bbb R}^n \longrightarrow {\Bbb R}$ of degree $m$, for any $\epsilon>0$, 
and any positive integer \break  $k \geq \gamma \epsilon^{-2} \delta^{-2} \ln (N+2)$, the 
inequality
$$ (1-\epsilon)^{m/2} \int_L f \ d\mu_k \leq 
 \left({k \over n}\right)^{m/2} \int_{{\Bbb R}^n} f \ d \mu_n \leq 
(1-\epsilon)^{-m/2} \int_L f \ d\mu_k$$
 holds with probability at least $2/3$ for a random $k$-dimensional subspace 
 $L \subset {\Bbb R}^n$.  
\endproclaim

Assuming that we can integrate efficiently over lower-dimensional subspaces (see Section 1.5 below), we get a randomized approximation algorithm for computing the integral of 
$f$ over ${\Bbb R}^n$. Namely, we sample a random $k$-dimensional subspace $L$,
compute the integral over $L$ and output the value of that integral multiplied by 
$(n/k)^{m/2}$.  To sample $L$ from the uniform 
 distribution on the Grassmannian $G_k({\Bbb R}^n)$, one can sample $k$ vectors
 $x_1, \ldots, x_k$ independently from the Gaussian 
distribution in ${\Bbb R}^n$ and let $L=\operatorname{span}\bigl\{ x_1, \ldots, x_k\bigr\}$.

  One ``anti-Monte Carlo'' feature of the algorithm is that the estimator is decidedly biased:
the expected value of the output  is essentially greater (by a factor of 
$(n/k)^{m/2}$) than the value we are trying to approximate. This is so because 
the distribution of the integral over a random subspace has a ``thick tail'': there are subspaces
which result in large integrals that significantly contribute to the integral 
over the whole space but such subspaces are very rare.

To increase the probability of obtaining the right approximation, one can use the standard 
approach of sampling several random subspaces and finding the median value of the outputs. 

One can observe that if $f$ is $(\delta, N)$-focused then $f^p$ is also $(\delta, N)$-focused 
for any positive integer $p$. This allows us to deduce that the maximum of $f$ over the 
unit sphere is well approximated by the scaled maximum of the restriction of $f$ onto the 
sphere in a lower-dimensional subspace.
\proclaim{(1.4) Corollary} There exists an absolute constant $\gamma>0$ with 
the following property.

For any $\delta >0$, for any positive integer $N$, for any  
$(\delta, N)$-focused polynomial \break $f: {\Bbb R}^n \longrightarrow {\Bbb R}$ of degree $m$, for any $\epsilon>0$, 
and any positive integer \break  $k \geq \gamma \epsilon^{-2} \delta^{-2} \ln (N+2)$, the 
inequality
$$(1-\epsilon)^{m/2}  \max_{x \in {\Bbb S}^{n-1} \cap L}  f(x) \leq 
\left({k \over n} \right)^{m/2} \max_{x \in {\Bbb S}^{n-1}} f(x)  \leq (1-\epsilon)^{-m/2}   
\max_{x \in {\Bbb S}^{n-1} \cap L}  f(x)$$
holds with probability at least $2/3$ for a random $k$-dimensional subspace 
$L \subset {\Bbb R}^n$.
\endproclaim
The problem of optimization of a polynomial on the unit sphere has attracted some attention 
recently, see \cite{F04} and \cite{K+04}.
Note that by restricting the polynomial onto a $k$-dimensional subspace we effectively reduce the number of variables to $k$ in the optimization problem. Using methods of computational 
algebraic geometry allows one to optimize a polynomial over the sphere in time 
exponential in the number of variables. Hence with $k=O(\log N)$, we obtain a quasi-polynomial algorithm of $m^{O(\log N)}$ complexity which approximates the maximum 
value of the polynomial on the sphere within a $(1-\epsilon)^{m/2}$ factor. If the degree $m$ of the polynomial
is fixed and $N$ is bounded by a polynomial in the number $n$ of variables, we get a polynomial 
time approximation algorithm.

\subhead (1.5) On the computational complexity\endsubhead
Let $f: {\Bbb R}^n \longrightarrow {\Bbb R}$ be a homogeneous polynomial 
of degree $m$ given by its ``black box'' which outputs the value of 
$f(x)$ for an input $x \in {\Bbb R}^n$. Then one can compute the 
monomial expansion 
$$f(x)=\sum_{\alpha} c_{\alpha} \xx^{\alpha} \quad \text{where} \quad \xx^{\alpha}=x_1^{\alpha_1} 
\ldots x_n^{\alpha_n} \quad \text{for} \quad \alpha=(\alpha_1, \ldots, \alpha_n)$$ 
in $O\left({n+m-1 \choose m}^3 \right)$ time through the standard procedure 
of interpolation, cf. also \cite{KY91} for the sparse version. If $L \subset {\Bbb R}^n$ is a
$k$-dimensional subspace, by choosing an orthonormal basis in $L$, we 
can identify $L$ with ${\Bbb R}^k$. Then the monomial expansion 
of the restriction $f_L$ can be computed in $O\left({k+m-1 \choose m}^3 \right)$ time. If $k$ is fixed, we get a polynomial time algorithm.
In we choose $k=O(\log N)$, the algorithms we obtain 
will be ``quasi-polynomial'', with the complexity of $m^{O(\log N)}$. 

Once a monomial expansion is obtained, it is easy to integrate polynomials
since there are explicit formulas to integrate monomials. Given 
a monomial $\xx^{\alpha}=x_1^{\alpha_1} \cdots x_n^{\alpha_n}$, the formula is
$$\int_{{\Bbb R}^n} \xx^{\alpha} \ d \mu_n = \cases  
\pi^{-n/2} \prod_{i=1}^n 2^{\alpha_i/2} \Gamma\left({\alpha_i+1 \over 2}\right) &\text{if all\ }
\alpha_i \ \text{are even} \\ 0 &\text{otherwise.} \endcases$$ 
\bigskip
In Section 2, we prove Theorem 1.3 and Corollary 1.4.
In Section 3, we consider some examples and applications, including 
the problem of approximating the {\it hafnian} of a positive matrix.
In Section 4, we consider the problem of integrating polynomials 
with respect to the complex Gaussian measure in ${\Bbb C}^n$. We prove a version of 
Theorem 1.3 in this case and show connections between efficient complex integration and
certain hard problems of combinatorial enumeration. 

\head 2. Proofs \endhead

One major ingredient of the proof of Theorems 1.3 is
the formula for the integral a product of linear forms.
\definition{(2.1) Definitions}
Let $m=2k$ be an even positive integer. A {\it perfect matching} $I$ of 
the set $\{1, \ldots, m\}$ is an unordered partition of
$\{1, \ldots, m\}$ into a union of $k$ unordered pairwise disjoint pairs 
$$I=\Bigl\{\{i_1, j_1\}, \{i_2, j_2\}, \ldots, \{i_k, j_k\} \Bigr\}.$$

Let $C=(c_{ij})$ be an $m \times m$ matrix, where $m=2k$ is an even integer.
The {\it hafnian} $\haf A$ of $A$ is defined by the formula
$$\haf C= \sum_I c_I,$$
where the sum is taken over all perfect matchings $I$ of the 
set $\{1, \ldots, m\}$ and $c_I$ 
is the product of all $c_{ij}$ for all pairs $\{i,j\} \in I$.
\enddefinition
The following result is known as  the
{\it Wick formula}, see, for example, \cite{Zv97}.
\proclaim{(2.2) Lemma} Let $m$ be a positive even integer and let 
$\ell_i: {\Bbb R}^n \longrightarrow {\Bbb R}$, $i=1, \ldots, m$,
be linear functions. Let 
$C=(c_{ij})$ be an $m \times m$ matrix defined by
$$c_{ij}=\int_{{\Bbb R}^n} \ell_i(x) \ell_j(x) \ d\mu_n.$$
Then
$$\int_{{\Bbb R}^n} \prod_{i=1}^m \ell_i(x) \ d \mu_n =
\haf C.$$
If $\ell_i$ is defined by $\ell_i(x)= \langle a_i, x \rangle$
for some $a_i \in {\Bbb R}^n$ then 
$c_{ij}=\langle a_i, a_j \rangle$.
\endproclaim
We also need a version of the  Johnson-Lindenstrauss ``flattenning''  Lemma, see, for example, \cite{Ve04}.
We present such a version below (with non-optimal constants), taken off Section V.7 of 
\cite{B02a}. 
\proclaim{(2.3) Lemma} Let $x \in {\Bbb R}^n$ be a vector and let 
$L \subset {\Bbb R}^n$ be a $k$-dimensional subspace 
chosen at random with respect to the Haar probability measure on 
the Grassmannian $G_k({\Bbb R}^n)$. Let $x'$ be the orthogonal projection
of $x$ onto $L$. Then, for any $0<\epsilon<1$, the probability that
$$(1-\epsilon) \|x\| \leq \sqrt{n \over k} \|x'\| \leq 
(1-\epsilon)^{-1} \|x\| $$
is at least $1-4\exp\{-\epsilon^2 k/4\}$.
\endproclaim

The following is a straightforward corollary. We establish it in a slightly larger generality than 
immediately needed, having in mind applications to complex integration in Section 4.
\proclaim{(2.4) Lemma} Let us choose $\delta>0$ and $\epsilon>0$. 
Suppose that $a_1, \ldots, a_N$ and $b_1, \ldots, b_N$ 
are vectors from ${\Bbb R}^n$ such that the cosine of the angle between
every pair $a_i$ and $b_j$ of vectors is at least $\delta>0$.

Let us choose a $\rho>0$ such that 
$$(1-\rho)^{-2} \leq 1+ {\delta \epsilon \over 3}$$ and an
integer 
$$k \geq \min \Bigl\{n, \quad 4\rho^{-2} \ln \Bigl(12 N^2 +24 N \Bigr)
\Bigr\}.$$
Let $L \subset {\Bbb R}^n$ be a $k$-dimensional subspace chosen 
at random with respect to the Haar probability measure on the Grassmannian
$G_k({\Bbb R}^n)$. Let $a_i', b_j'$ be the orthogonal projection 
of $a_i$, $b_j$ onto $L$. Then with probability at least $2/3$    
$$(1-\epsilon) \langle a_i, b_j \rangle \leq {n \over k} 
\langle a_i', b_j' \rangle \leq (1-\epsilon)^{-1} \langle a_i, b_j \rangle$$
for all pairs $(i,j)$.
\endproclaim
\demo{Proof} Scaling, if necessary, we may assume that 
$\|a_i\|=\|b_j\|=1$ for all $i$ and $j$, so 
$\langle a_i, b_j \rangle \geq \delta$ for all $i,j$. 
We have
$$\langle a_i, b_j \rangle=
{\|a_i + b_j\|^2 - \|a_i\|^2 - \|b_j\|^2 \over 2} \quad 
\text{and} \quad \langle a_i', b_j' \rangle=
{\|a_i'+b_j'\|^2 -\|a_i'\|^2 -\|b_j'\|^2 \over 2}.$$
We note that
$$(1-\rho)^{-2} \leq 1+ {\delta \epsilon \over 3} \quad \text{and} \quad 
(1-\rho)^2 \geq 1-{\delta \epsilon \over 3}.$$
Since there are altogether $N^2 +2N$ vectors $a_i, b_j$, and $a_i+b_j$,
by Lemma 2.3, for a random $k$-dimensional subspace $L$, with 
probability at least $2/3$, we get
$$\|a_i+b_j\|^2 (1-\rho)^2 \leq {n \over k} \|a_i' +b_j'\|^2 
\leq (1-\rho)^{-2} \|a_i + b_j\|^2$$
and, similarly,
$$\split &\|a_i\|^2 (1-\rho)^2 \leq {n \over k} \|a_i'\|^2 
\leq (1-\rho)^{-2} \|a_i\|^2 \quad \text{and} \\
&|b_i\|^2 (1-\rho)^2 \leq {n \over k} \|b_i'\|^2 
\leq (1-\rho)^{-2} \|b_i\|^2 \endsplit$$
for all pairs $i,j$.
Since $\|a_i\|=\|b_j\|=1$ and $\| a_i + b_j\| \leq 2$, we get
$$\|a_i +b_j \|^2 -{4 \delta \epsilon \over 3}
 \leq {n \over k} \|a_i' +b_j'\|^2 
\leq \|a_i+b_j\|^2 +{4 \delta \epsilon \over 3}$$
and, similarly,
$$\split &\|a_i\|^2 -{\delta \epsilon \over 3} \leq {n \over k} \|a_i'\|^2 
\leq  \|a_i\|^2 +{\delta \epsilon \over 3} \quad \text{and} \\
&|b_i\|^2 -{\delta \epsilon \over 3} \leq {n \over k} \|b_i'\|^2 
\leq \|b_i\|^2 +{\delta \epsilon \over 3}. \endsplit$$
Therefore,
$$\langle a_i, b_j \rangle -\delta \epsilon \leq {n \over k} 
\langle a_i', b_j' \rangle  \leq \langle a_i, b_j \rangle + \delta \epsilon.$$
Since $\langle a_i, b_j \rangle \geq \delta$, the 
proof follows.
{\hfill \hfill \hfill} \qed
\enddemo
\proclaim{(2.5) Corollary} There exists an absolute constant $\gamma>0$ 
with the following property.

Let  $\delta>0$ and $\epsilon>0$ be numbers, let $N$ be a positive integer, and
let   $a_1, \ldots, a_N$ and $b_1, \ldots, b_N$ be vectors from 
${\Bbb R}^n$ such that the cosine of the angle between 
every pair $a_i, b_j$ of vectors is at least $\delta$. Let $k$ be a positive integer such that
$$k \geq  \gamma \delta^{-2} \epsilon^{-2} \ln (N+2) $$
and let $L \subset {\Bbb R}^n$ be a $k$-dimensional subspace chosen at random
with respect to the Haar probability measure in the Grassmannian $G_k({\Bbb R}^n)$.
Let $a_i', b_j'$ be the orthogonal projections of $a_i, b_j$ onto $L$. 
Then, with probability
at least $2/3$, we have
$$(1-\epsilon) \langle a_i', b_j' \rangle \leq {k \over n}
\langle a_i, b_j \rangle \leq (1-\epsilon)^{-1} \langle a_i', b_j' \rangle $$
for all pairs $a_i, b_j$.    
\endproclaim
The proof follows by Lemma 2.4.
\medskip
Now we are ready to prove Theorem 1.3.
\medskip
\demo{Proof of Theorem 1.3}
 We can write 
 $$f(x)=\sum_I \alpha_I  \prod_{i \in I} \langle c_i, x \rangle,$$
 where the cosine of the angle between every pair of vectors $c_i$ and $c_j$ is at least 
 $\delta$, $I$ ranges over subsets $I \subset \{1, \ldots, N\}$ of cardinality $m$, and 
 $\alpha_I \geq 0$. 
 For every $I$, let us consider the $m \times m$ matrix $C_I$ whose entries $c_{ij}$ are defined 
 by $c_{ij}=\langle c_i, c_j \rangle$.
 Then, by Lemma 2.2,
 $$\int_{{\Bbb R}^n} f(x) \ d\mu_n=\sum_I \alpha_I \haf C_I.$$
 Let $L \subset {\Bbb R}^n$ be a $k$-dimensional subspace. 
 Then the restriction $f_L$ of $f$ onto $L$ can be 
 written as 
 $$f_L(x)=\sum_I \alpha_I  \prod_{i \in I} \langle c_i', x \rangle,$$
 where $c_i'$ are the orthogonal projections of $c_i$ onto $L$.
 Therefore,
 $$\int_L f(x) \ d \mu_k =\sum_I \alpha_I \haf C_I', $$
 where the entries $c_{ij}'$ of $C_I'$ are defined by $c_{ij}'=\langle c_i', c_j' \rangle$.
 Since the hafnian of an $m \times m$ matrix is a non-negative homogeneous polynomial of 
 degree $m/2$ in the entries of the matrix, the proof follows by Corollary 2.5 where 
 we take $a_i=b_i=c_i$.
 {\hfill \hfill \hfill} \qed
\enddemo

\demo{Proof of Corollary 1.4} First, we claim that 
$$\max_{x \in {\Bbb S}^{n-1}} f(x) = \max_{x \in {\Bbb S}^{n-1}} |f(x)|.$$
If the degree $m$ of $f$ is odd, this is immediate. If $m$ is even, let us consider 
the polynomial $f^p$ for some odd $p$. Since
$$f(x)=\sum_I \alpha_I \prod_{i \in I} \langle c_i, x \rangle \quad \text{where} \quad \alpha_I \geq 0,$$
the polynomial $f^p$ is also represented as a non-negative linear combination of 
products of $\langle c_i, x \rangle$, where the cosine of the angle between every pair 
$c_i, c_j$ of vectors is at least $\delta$. It follows from the proof of Theorem 1.3 above that 
$$\int_{{\Bbb S}^{n-1}} f^p \ d x >0 \quad \text{for any} \quad p.$$
from which we conclude that the maximum value of $f$ and the maximum absolute value of $f$ on the sphere ${\Bbb S}^{n-1}$
must coincide.

Next, as in the proof of Theorem 1.3,  we observe that if $L \subset {\Bbb R}^n$ is a 
$k$-dimensional subspace such that for the orthogonal projections
$c_1', \ldots, c_N'$ of $c_1, \ldots, c_N$ onto $L$ we have
$$(1-\epsilon) \langle c_i', c_j' \rangle \leq {k \over n} \langle c_i, c_j \rangle \leq
(1-\epsilon)^{-1} \langle c_i', c_j' \rangle \quad \text{for all pairs} \quad i,j $$    
Then 
$$(1-\epsilon)^{mp/2} \int_L f^p \ d\mu_k \leq \left(k \over n\right)^{mp/2} \int_{{\Bbb R}^n} 
f^p \ d\mu_n \leq (1-\epsilon)^{-mp/2} \int_L f^p \ d\mu_k $$
for all $p$. 
In particular, if the degree $m$ of $f$ is even, 
$$\int_{{\Bbb S}^{n-1} \cap L} f^p \ dx >0 \quad \text{for all} \quad 
p.$$
Therefore, 
$$\max_{x \in {\Bbb S}^{n-1} \cap L} f(x) =\max_{x \in {\Bbb S}^{n-1} \cap L} |f(x)|.$$
The proof now follows from the identities
$$\split &\lim_{p \longrightarrow +\infty} \left(\int_{{\Bbb S}^{n-1}} f^{2p}(x) \ dx\right)^{1/2p}=\max_{x \in {\Bbb S}^{n-1}} |f(x)|=\max_{x \in {\Bbb S}^{n-1}} f(x), \\ 
 &\lim_{p \longrightarrow +\infty} \left(\int_{{\Bbb S}^{n-1} \cap L} f^{2p}(x) \ dx\right)^{1/2p}=\max_{x \in {\Bbb S}^{n-1} \cap L} |f(x)|=\max_{x \in {\Bbb S}^{n-1} \cap L} f(x), \\
& \int_{{\Bbb S}^{n-1}} f^{2p}(x) \ dx =
{\Gamma(n/2) \over 2^{mp} \Gamma(n/2+mp)} \int_{{\Bbb R}^n} f^{2p} \ d\mu_n, \quad \text{and} \\
& \int_{{\Bbb S}^{n-1} \cap L} f^{2p}(x) \ dx =
{\Gamma(k/2) \over 2^{mp} \Gamma(k/2+mp)} \int_{L} f^{2p} \ d\mu_k. \endsplit $$ 
\enddemo
{\hfill \hfill \hfill} \qed

\head 3. Examples and an Application \endhead

Some natural examples of sets of vectors $c_1, \ldots, c_N \in {\Bbb R}^n$ with the 
property that for every $(i,j)$, the cosine of the angle between $c_i$ and $c_j$ is at least 
$\delta>0$ are as follows.
\example{(3.1) Examples}

(3.1.1)  Let $c_1, \ldots, c_N \in {\Bbb R}^n$ 
be vectors with positive coordinates such that the ratio of the smallest/largest coordinate 
for each vector $c_i$ is at least $\sqrt{\delta}$. It is easy to show that the cosine of the angle 
between $c_i$ and $c_j$ is at least $\delta$ for each pair $(i,j)$.
\medskip
(3.1.2) Suppose that $n=k(k+1)/2$ and let us identify ${\Bbb R}^n$ with the space of 
$k \times k$ symmetric matrices with the scalar product $\langle a, b \rangle=\text{trace}(ab)$.
Let $c_1, \ldots, c_N$ be positive definite matrices such that the ratio of the smallest/largest
eigenvalue for each matrix $c_i$ is at least $\sqrt{\delta}$. It is easy to show that the
cosine of the angle between $c_i$ and $c_j$ is at least $\delta$ for each pair $(i,j)$. 
\medskip
Other examples can be obtained by sampling $c_1,  \ldots, c_N$ at random from some 
biased distribution in ${\Bbb R}^n$ (a distribution with a non-zero expectation).
\endexample 

Whenever we have a polynomial 
$$f(x)=\sum \Sb I \subset \{1, \ldots, N \} \\ |I|=m \endSb  \alpha_I \prod_{i \in I} 
\langle c_i, x \rangle \quad \text{where} \quad \alpha_I \geq 0$$
and vectors $c_i$ as in (3.1.1)-(3.1.2), integration (optimization) of such a polynomial over 
the unit sphere ${\Bbb S}^{n-1}$ reduces to integration (optimization) over a random lower-dimensional
subspace $L$. If we want to achieve a $(1-\epsilon)^m$ factor of approximation, the 
dimension $k$ of the subspace is only logarithmic in $N$, so that as long as $N$ is bounded 
by a polynomial in $n$, we achieve an exponential reduction in the number of variables.  

Finally, we consider the problem of computing (approximating) the hafnian of a given positive 
matrix. This problem is of interests in combinatorics and statistical physics and generalizes the 
problem of computing the permanent, see Section 8.2 of \cite{Mi78}. Unlike in the case of the 
permanent, where a 
polynomial time approximation algorithm has been recently obtained \cite{J+04}, much less
is known about computing hafnians.

\subhead (3.2) Computing the hafnian of a positive matrix \endsubhead
Let $C=(c_{ij})$ be an $m \times m$ positive symmetric matrix, where 
$m=2k$ is even. Recall (see Definition 2.1) that the 
hafnian of $C$ is the polynomial
$$\haf C=\sum_I  c_I,$$
where the sum is taken over all perfect matchings 
$I=\Bigl\{ \{i_1, j_1\}, \ldots, \{i_k, j_k\} \Bigr\}$ of the set
$\{1, \ldots, m\}$ and $c_I$ is the product of $c_{ij}$ for 
$\{i, j\} \in I$. 

Suppose that $C$ is positive semidefinite. Then $C$ is the Gram matrix 
of a set of vectors, so $c_{ij}=\langle c_i, c_j \rangle$ for 
some vectors $c_1, \ldots, c_m \in {\Bbb R}^m$ and such a representation 
can be computed efficiently (in polynomial time). Using the Wick formula 
(Lemma 2.2), we can write 
$$\haf C=\int_{{\Bbb R}^m} \prod_{i=1}^m \langle c_i, x \rangle \ d\mu_m.$$
Suppose that for each pair $c_i, c_j$ of vectors the cosine of the angle
between $c_i$ and $c_j$ is at least $\delta$, which means that 
$c_{ij} \geq \delta \sqrt{c_{ii} c_{jj}}$ for every pair $i,j$. Then, 
by Theorem 1.3, to approximate $\haf C$ within a factor of 
$(1-\epsilon)^{m/2}$, we can replace the integral by the integral over 
a random $k$-dimensional subspace $L \subset {\Bbb R}^m$ with 
$k=O\bigl(\epsilon^{-2} \delta^{-2} \ln (m+2)\bigr)$. If $\epsilon$ and $\delta$ 
are fixed in advance, we get a quasi-polynomial algorithm of 
$m^{O(\ln m)}$ complexity.

One can extend the above argument as follows. We observe that $\haf C$ 
does not depend at all on the diagonal entries of $C$, so we are free 
to change the diagonal entries of $C$ to ensure that the above conditions
are satisfied. If we put sufficiently large numbers on the diagonal of 
$C$, we can make sure that $C$ is positive definite, so 
$c_{ij}=\langle c_i, c_j \rangle$ for some vectors 
$c_1, \ldots, c_m \in {\Bbb R}^m$. The goal is to make the cosine of 
the angle between every pair $c_i, c_j$ of vectors as large as possible.
Suppose that $c_{ii}=0$ for all $i$ and let $-\lambda$ be the minimum
eigenvalue of $C$. Then $C+\lambda I$ is a positive semidefinite matrix 
and the cosine of the angle between $c_i$ and $c_j$ is $c_{ij}/\lambda$.
Thus as long as the absolute value $\lambda$ of negative 
eigenvalues of $C$ is sufficiently small, we get an efficient algorithm 
to approximate $\haf C$.

\head 4. Complex Integration \endhead

Let $f, g: {\Bbb R}^n \longrightarrow {\Bbb R}$ be real $n$-variate homogeneous polynomials. Let us identify ${\Bbb R}^n \oplus {\Bbb R}^n = {\Bbb C}^n$ via 
$x+iy=z$ and let $\nu_n$ be the Gaussian measure on ${\Bbb C}^n$ with the density 
$$\pi^{-n} e^{-\|z\|^2}, \quad \text{where} \quad \|z\|^2=\|x\|^2 + \|y\|^2 \quad \text{for} \quad 
z=x+iy.$$
We recall that $\overline{z}=x-iy$ is the complex conjugate of $z=x+iy$.

Let us define the scalar product on the space of polynomials
$$\langle f, g \rangle =\int_{{\Bbb C}^n} f(z) \overline{g(z)} \ d \nu_n$$
(although we use the same notation for the standard scalar product on ${\Bbb R}^n$, we 
hope no confusion will result since the domains are drastically different).
One can easily check that the monomials 
$$ \xx^{\alpha}=x_1^{\alpha_1} \ldots x_n^{\alpha_n} \quad \text{for} \quad 
\alpha=(\alpha_1, \ldots, \alpha_n), 
\quad \text{where} \quad \alpha_i \geq 0 \quad
\text{for} \quad i=1, \ldots, n.$$
are orthogonal under the scalar product, though not orthonormal:
$$\langle \xx^{\alpha}, \xx^{\beta} \rangle=\cases \alpha_1! \ldots \alpha_n! &\text{if \ }
\alpha=\beta=(\alpha_1, \ldots, \alpha_n) \\ 0 &\text{if\ } \alpha \ne \beta. \endcases$$
Therefore, if 
$$f=\sum_{\alpha \in F} a_{\alpha} \xx^{\alpha} \quad \text{and} \quad 
g=\sum_{\alpha \in G} b_{\alpha} \xx^{\alpha}$$
are the monomial expansions of $f$ and $g$, we have 
$$\langle f, g \rangle =\sum_{\alpha \in F \cap G} a_{\alpha} b_{\alpha} \alpha_1 ! \cdots \alpha_n!.$$  
It follows from the integral representation that the scalar product is invariant under the action of
the orthogonal group: if $U$ is an orthogonal transformation of ${\Bbb R}^n$ and 
polynomials $f_1$, $g_1$ are defined by $f_1(x)=f(Ux)$ and $g_1(x)=g(Ux)$, then
$\langle f_1, g_1 \rangle=\langle f, g \rangle$. 

Various problems of combinatorial counting reduce to computing the scalar
products of two polynomials.
\example{(4.1) Example} Let $a_1, \ldots, a_N$ and $b$ be some non-negative
integer $n$-vectors. Let $M$ be a positive integer. We define 
$$f(x)=\prod_{i=1}^N \left(\sum_{k=0}^M \xx^{k a_i} \right) \quad 
\text{and} \quad g(x)= \xx^b.$$
Then the monomial expansion of $f$ contains all monomials $\xx^{a}$, where 
$a$ is a linear combination of $a_1, \ldots, a_N$ with positive integer  
coefficients not exceeding $M$. Furthermore, if 
$b=(\beta_1, \ldots, \beta_n)$, then 
$\langle f, g \rangle$ is the number of non-negative integer solutions
$(k_1, \ldots, k_N)$, $0 \leq k_i \leq M$, to the equation
$$k_1 a_1 + \ldots + k_N a_N = b$$
times $\beta_1 ! \ldots \beta_n!$. 
The number of such solutions $(k_1, \ldots, k_N)$ as a function of $b$ is often called
the {\it vector partition function}, cf. \cite{BV97}. Computing the vector partition function is generally  
as hard as counting integer points in a polytope.   
\endexample 

\definition{(4.2) Definition} Let us fix a number $0< \delta \leq 1$ and a positive integer $N$. We say 
that a pair of homogeneous polynomials $f, g: {\Bbb R}^n \longrightarrow {\Bbb R}$ 
of degree $m$ is $(\delta, N)$-{\it focused} if 
there exist $N$ non-zero vectors $a_1, \ldots, a_N \in {\Bbb R}^n$ and $N$ non-zero
vectors $b_1, \ldots, b_N \in {\Bbb R}^n$ such that  
\medskip
$\bullet$ for every pair $(i,j)$ the cosine of the angle between $a_i$ and $b_j$ is at least 
$\delta$;
\medskip
$\bullet$ the polynomial $f$ can be written as a non-negative linear combination 
$$f(x)=\sum_I \alpha_I \prod_{i \in I} \langle a_i, x \rangle,$$
while the polynomial $g$ can be written as a non-negative linear combination 
$$g(x)=\sum_I \beta_I \prod_{j \in J} \langle b_j, x \rangle,$$
where the sum is taken over subsets $I, J \subset \{1, \ldots, m\}$ of cardinality 
$|I|=|J|=m$ and $\alpha_I, \beta_J \geq 0$.   
\enddefinition

We prove that the value of the scalar product of a well-focused pair
of polynomials can be well-approximated from the scalar product of the restriction of the polynomials onto a random lower-dimensional subspace.

For a $k$-dimensional subspace $L \subset {\Bbb R}^n$, let us consider
its complexification $L_{\Bbb C}=L \oplus i L \subset {\Bbb C}^n$. Let 
$\nu_k$ be the Gaussian measure in $L_{\Bbb C}$ with the density
$\pi^{-k} \exp\{-\|z\|^2\}$ for $z \in L_{\Bbb C}$. We pick a $k$-dimensional subspace 
$L \subset {\Bbb R}^n$ at random with respect to the Haar probability 
measure on the Grassmannian $G_k({\Bbb R}^n)$ and consider the restrictions $f_L$ and $g_L$ onto
$L$ and  the integral
$$\langle f_L, g_L \rangle=
\int_{L_{\Bbb C}} f(z) \overline{g(z)} \ d \nu_k.$$ 
We claim that as long as $k \sim \log N$, the properly scaled integral over $L_{\Bbb C}$ approximates
the  integral over ${\Bbb C}^n$ within a factor of $(1-\epsilon)^m$.

\proclaim{(4.3) Theorem} There exists an absolute constant $\gamma>0$ with 
the following property.

For every $\delta>0$, for any positive integer $N$, for any 
$(\delta, N)$-focused pair of polynomials 
$f,g: {\Bbb R}^n \longrightarrow {\Bbb R}$ of degree $m$, for any 
$\epsilon >0$ and any positive integer 
$k \geq \gamma \epsilon^{-2} \delta^{-2} \ln (N+2)$, 
the inequality
$$(1-\epsilon)^m  
\langle f_L, g_L \rangle \leq \left({k \over n}\right)^m \langle f, g \rangle 
\leq (1-\epsilon)^{-m}  
\langle f_L, g_L \rangle$$
holds with probability at least $2/3$ for a random $k$-dimensional
subspace $L \subset {\Bbb R}^n$.  
\endproclaim

The proof is very similar to that of Theorem 1.3. The only difference
is that we need the complex version of the Wick formula.

\definition{(4.4) Definitions}
Let $m$ be a positive integer. A {\it permutation} of the set \break
$\{1, \ldots, m\}$ is a bijection  
$\sigma: \{1, \ldots, m\} \longrightarrow \{1, \ldots, m\}$.

Let $C=(c_{ij})$ be an $m \times m$ matrix.
The {\it permanent} $\per C$ of $C$ is defined by the 
formula
$$\per C=\sum_{\sigma} \prod_{i=1}^m c_{i \sigma(i)},$$
where the sum is taken over all permutations of the set 
$\{1, \ldots, m\}$.
\enddefinition
Here is the complex version of the Wick formula. Since the author was unable to locate 
it in the literature, a proof is given here.  
\proclaim{(4.5) Lemma} Let $m$ be a positive integer and let 
$f_i, g_i: {\Bbb R}^n \longrightarrow {\Bbb R}$ be linear functions.
Let $C=(c_{ij})$ be an $m \times m$ matrix defined by 
$$c_{ij}=\int_{{\Bbb C}^n} f_i(z) \overline{g_j(z)} \ d \nu_n.$$
Then
$$ \int_{{\Bbb C}^n} \prod_{i=1}^n f_i(z) \overline{g_i(z)} \ d \nu_n=
\per C.$$  
If $f_i$ is defined by $f_i(x)=\langle a_i, x \rangle$ and
$g_j$ is defined by $g_i(x)=\langle b_j, x \rangle$ for some
$a_i, b_j \in {\Bbb R}^n$ then $c_{ij}=\langle a_i, b_j \rangle$.
\endproclaim   
\demo{Proof}
Given vectors $a_1, \ldots, a_m$ and $b_1, \ldots, b_m$, let 
$$p(x)=\prod_{i=1}^m \langle a_i, x \rangle \quad \text{and} \quad 
q(x)=\prod_{j=1}^m \langle b_j, x \rangle.$$
Our goal is to prove that 
$$\langle p, q \rangle =\per C \quad \text{where} \quad 
c_{ij}=\langle a_i, b_j \rangle.$$

First, we check the identity in the special case when
$a_1= \ldots = a_m =e_1$, the first basis vector, and 
$b_1 =\ldots = b_m =b=(\beta_1, \ldots, \beta_n)$ is an arbitrary vector.
In this case, 
$p(x)=x_1^m$ and $q(x)=(\beta_1 x_1 + \ldots + \beta_n x_n)^m$,
so we have $\langle p, q \rangle=\beta_1^m m!$.
 On the other hand, $c_{ij}=\beta_1$ for all $i$ and $j$,
so $\per C=m! \beta_1^m$ as well.

Next, we check the identity when $a_1, \ldots, a_m=a$ and $b_1, \ldots, b_m=b$, where 
$a$ and $b$ are arbitrary vectors. Applying scaling, if necessary, we can assume that 
$\|a\|=1$. Since an orthogonal transformation of ${\Bbb R}^n$ does not change 
either $\langle p, q \rangle$ or $C$, this case reduces to the previous one.

Now we consider the general case. We observe that both quantities 
$\langle p, q \rangle$ and $\per C$ are multilinear and symmetric in 
$a_1, \ldots, a_m$  and multilinear and symmetric in $b_1, \ldots, b_m$, so we obtain the general 
case by polarization.
For variables $\lambda=(\lambda_1, \ldots, \lambda_m)$ and $\mu=(\mu_1, \ldots, \mu_m)$
we introduce vectors $a_{\lambda}=\lambda_1 a_1 + \ldots + \lambda_m a_m$ 
and $b_{\mu}=\mu_1 b_1 + \ldots + \mu_m b_m$. 
If $F(a_1, \ldots, a_m; b_1, \ldots, b_m)$ is any polynomial multilinear and symmetric in 
$a_1, \ldots, a_m$ and multilinear and symmetric in $b_1, \ldots, b_m$, then 
$(m!)^2 F(a_1, \ldots, a_m; b_1, \ldots, b_m)$ is equal to the coefficient of the product 
$\lambda_1 \cdots \lambda_m \mu_1 \cdots \mu_m$ in the expansion of 
 $F(a_{\lambda}, \ldots, a_{\lambda}; b_{\mu}, \ldots, b_{\mu})$ as a polynomial 
 in $\lambda_1, \ldots, \lambda_m, \mu_1, \ldots, \mu_m$. Since if two such polynomials
 $F$ and $G$ agree on all $(2m)$-tuples $(a, \ldots, a; b, \ldots, b)$, they agree everywhere.
 Letting $F=\langle p, q \rangle$ and $G=\per C$, we complete the proof.
 {\hfill \hfill \hfill} \qed
\enddemo

Now the proof of Theorem 4.5 follows the proof of Theorem 1.3.

\head References \endhead

\widestnumber\key{MMMM}

\ref\key{B02a}
\by A. Barvinok 
\book A Course in Convexity 
\bookinfo Graduate Studies in Mathematics
\vol  54 
\publ American Mathematical Society
\publaddr Providence, RI
\yr  2002
\endref

\ref\key{B02b}
\by A. Barvinok
\paper Estimating $L\sp \infty$ norms by $L\sp {2k}$ norms for functions on orbits
\jour  Found. Comput. Math. 2 
\yr 2002
\pages 393--412
\endref

\ref\key{BV97}
\by M. Brion and M. Vergne
\paper Residue formulae, vector partition functions and lattice points in rational polytopes
\jour  J. Amer. Math. Soc.
\vol 10 
\yr 1997
\pages  797--833
\endref

\ref\key{Fa04} 
\by ÊL. Faybusovich
\paper Global optimization of homogeneous polynomials on the simplex and on the sphere
\inbook  Frontiers in global optimization
\pages 109--121
\bookinfo Nonconvex Optim. Appl.
\vol  74
\publ  Kluwer Acad. Publ.
\publaddr Boston, MA
\yr  2004
\endref

\ref\key{J+04}
\by M. Jerrum, A. Sinclair, and E. Vigoda
\paper A polynomial-time approximation algorithm for the permanent of a matrix with 
non-negative entries
\jour Journal of the ACM
\vol 51
\yr 2004
\pages 671--697
\endref

\ref\key{KY91}
\by E. Kaltofen and L. Yagati
\paper Improved sparse multivariate polynomial interpolation algorithms
\inbook  Symbolic and algebraic computation (Rome, 1988)
\pages 467--474
\paperinfo Lecture Notes in Comput. Sci.
\vol 358
\publ Springer
\publaddr Berlin
\yr 1989
\endref

\ref\key{K+04}
\by E. De Klerk, M. Laurent, and P. Parrilo
\paper A PTAS for the minimization of polynomials of fixed degree over the simplex
\paperinfo preprint
\yr 2004
\endref

\ref\key{Mi78}
\by H. Minc
\book Permanents 
\bookinfo Encyclopedia of Mathematics and its Applications
\vol  6  
\publ Addison-Wesley Publishing Co.
\publaddr Reading, Mass.
\yr 1978
\endref

\ref\key{MS86} 
\by V.D. Milman and  G. Schechtman
\book Asymptotic Theory of Finite-\ Dimensional Normed Spaces. With an Appendix by M. Gromov
\bookinfo Lecture Notes in Mathematics
\vol 1200
\publ Springer-Verlag
\publaddr  Berlin
\yr  1986
\endref

\ref\key{Re92}
\by J. Renegar
\paper  On the computational complexity and geometry of the first-order theory of the reals
\jour J. Symbolic Comput. 
\vol 13 
\yr 1992
\pages 255--352
\paperinfo I. Introduction. Preliminaries. The geometry of semi-algebraic sets. The decision problem for the existential theory of the reals (255--299), II. The general decision problem. Preliminaries for quantifier elimination (301--327), III. On the computational complexity and geometry of the first-order theory of the reals. III. Quantifier elimination (329--352)
\endref

\ref\key{Ve04}
\by S.S. Vempala 
\book The Random Projection Method 
\bookinfo DIMACS Series in Discrete Mathematics and Theoretical Computer Science
\vol 65 
\publ American Mathematical Society
\publaddr  Providence, RI
\yr 2004
\endref

\ref\key{Zv97}
\by A. Zvonkin
\paper Matrix integrals and map enumeration: an accessible introduction
\paperinfo Combinatorics and physics (Marseilles, 1995) 
\jour Math. Comput. Modelling 
\vol 26 
\yr 1997
\pages  281--304
\endref

\enddocument

\end